\documentclass[12pt]{article}
\usepackage{amsmath,amssymb}
\usepackage{geometry}
\usepackage{amsmath, amsfonts,amsthm,amssymb,enumerate,amscd}

\usepackage{times}

\usepackage[active]{srcltx}
\usepackage[utf8]{inputenc}
\usepackage{algorithm}
\usepackage{stmaryrd}
\usepackage{array}
\usepackage{algorithmicx}
\usepackage{algpseudocode}
\usepackage{graphicx}
\usepackage{letltxmacro}

\usepackage{color}
\usepackage[table,xcdraw,svgnames,dvipsnames]{xcolor}
\definecolor{light-salmon}{RGB}{255,140,120}
\usepackage[colorlinks]{hyperref}
\hypersetup{
	allbordercolors=.,
	citecolor=DarkGreen,
	linkcolor=DarkRed,
	urlcolor=NavyBlue
}

\graphicspath{{./pics/}}

\usepackage{enumitem}

\numberwithin{equation}{section}
\theoremstyle{plain}
\newtheorem{thm}{Theorem}

\newtheorem{cor}[thm]{Corollary}
\newtheorem{prop}[thm]{Proposition}
\newtheorem{deff}[thm]{Definition}
\theoremstyle{definition}
\newtheorem{rem}[thm]{Remark}

\newcommand{\Area}{\operatorname{Area}}

\renewcommand{\Bbb}{\mathbb}
\newcommand{\bo}[1]{{\bf #1}}

\makeatletter
\DeclareFontFamily{U}{tipa}{}
\DeclareFontShape{U}{tipa}{m}{n}{<->tipa10}{}
\newcommand{\arc@char}{{\usefont{U}{tipa}{m}{n}\symbol{62}}}%

\newcommand{\arc}[1]{\mathpalette\arc@arc{#1}}

\newcommand{\arc@arc}[2]{%
	\sbox0{$\m@th#1#2$}%
	\vbox{
		\hbox{\resizebox{\wd0}{\height}{\arc@char}}
		\nointerlineskip
		\box0
	}%
}
\makeatother

\title{New variational arguments regarding the Blaschke-Lebesgue theorem}

\author{Beniamin Bogosel}

\begin{document}
	
\maketitle

\begin{abstract}
The sensitivity of the areas of Reuleaux polygons and disk polygons is computed with respect to vertex perturbations.  Computations are completed for both constrained and Lagrangian formulations and they imply that the only critical Reuleaux polygons for the area functional are the regular ones. As a consequence, new variational proofs for the Blaschke-Lebesgue and Firey-Sallee theorems are found.
\end{abstract}



\section{Introduction}

\begin{figure}
	\centering
	\includegraphics[width=0.25\textwidth]{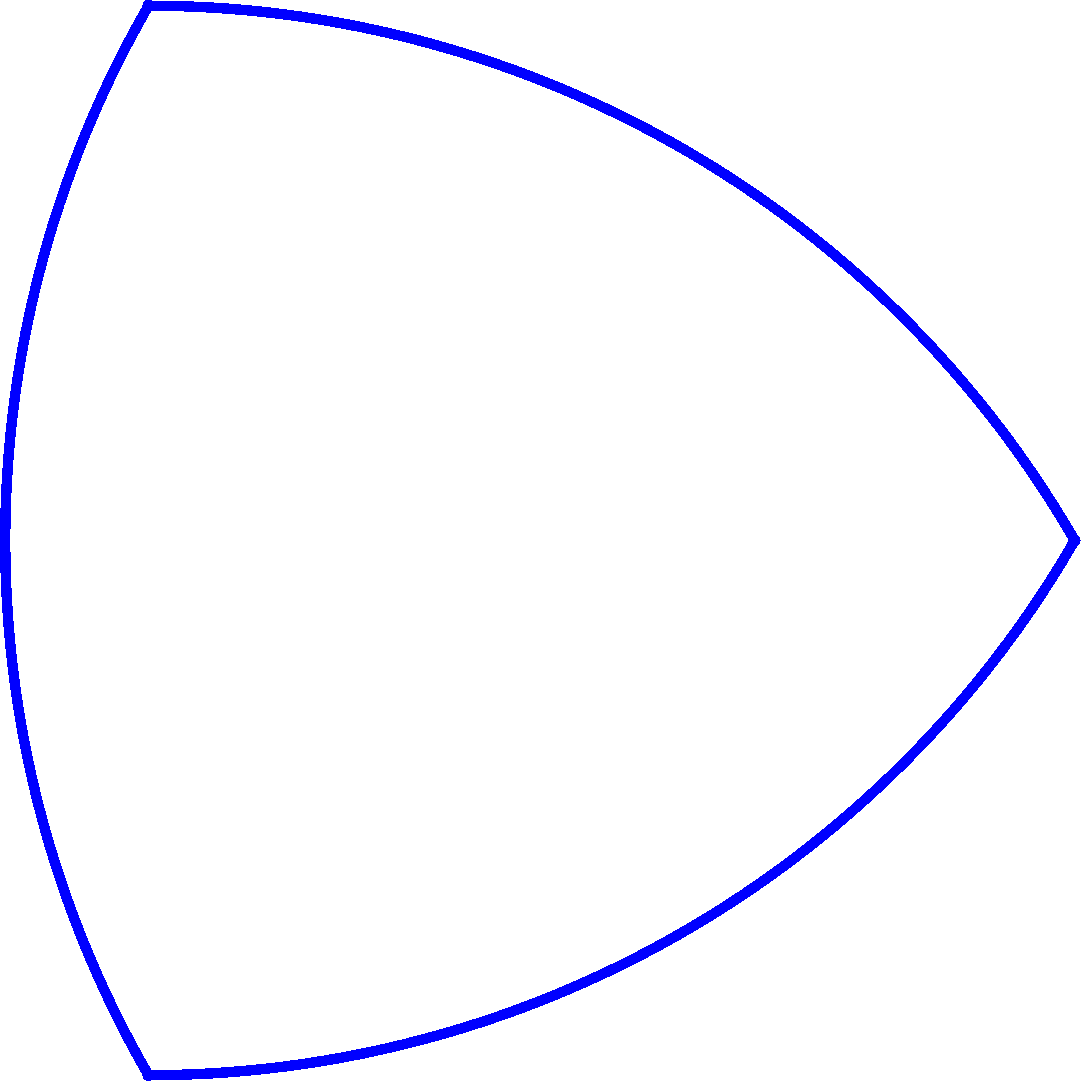}\quad
	\includegraphics[width=0.25\textwidth]{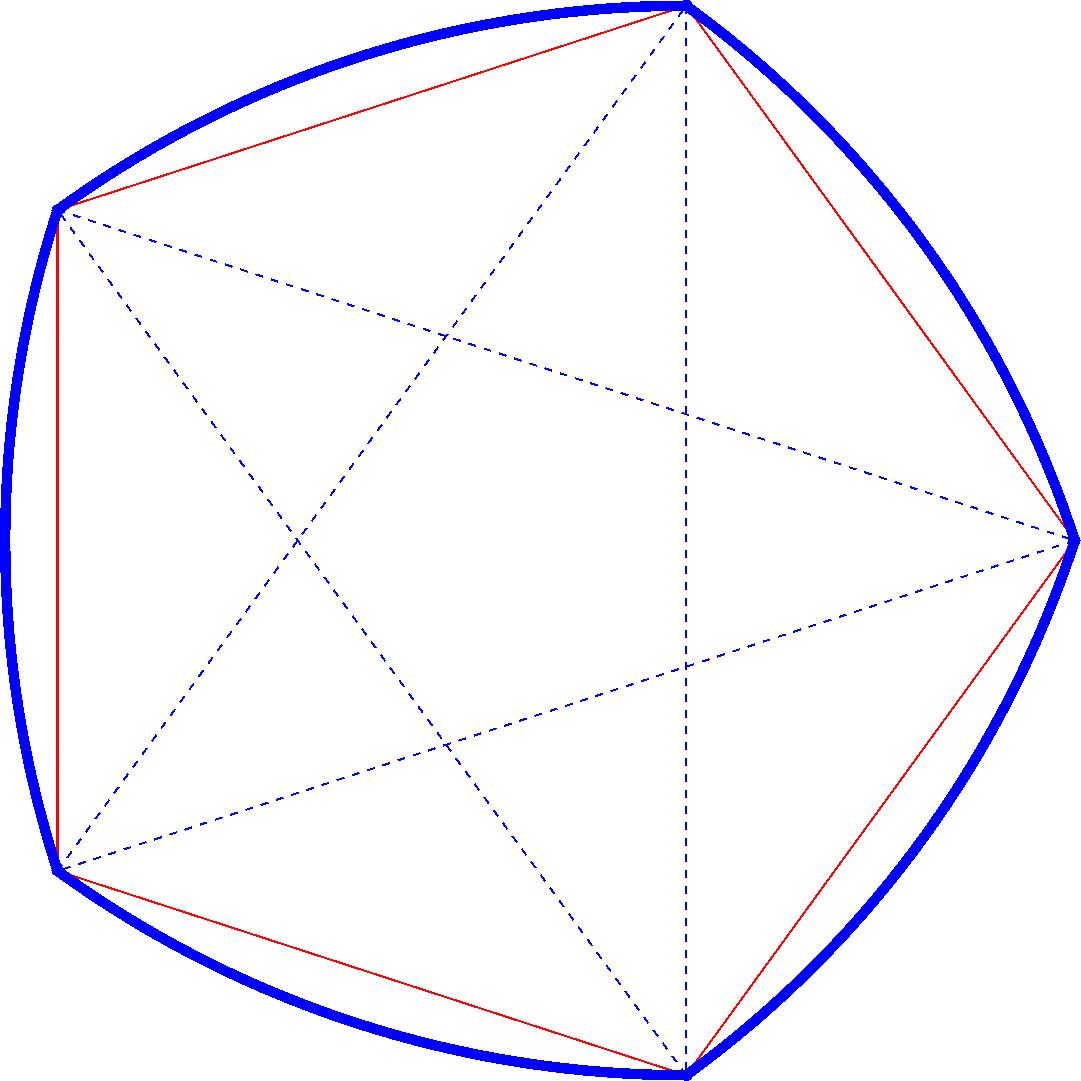}\quad
	\includegraphics[width=0.25\textwidth]{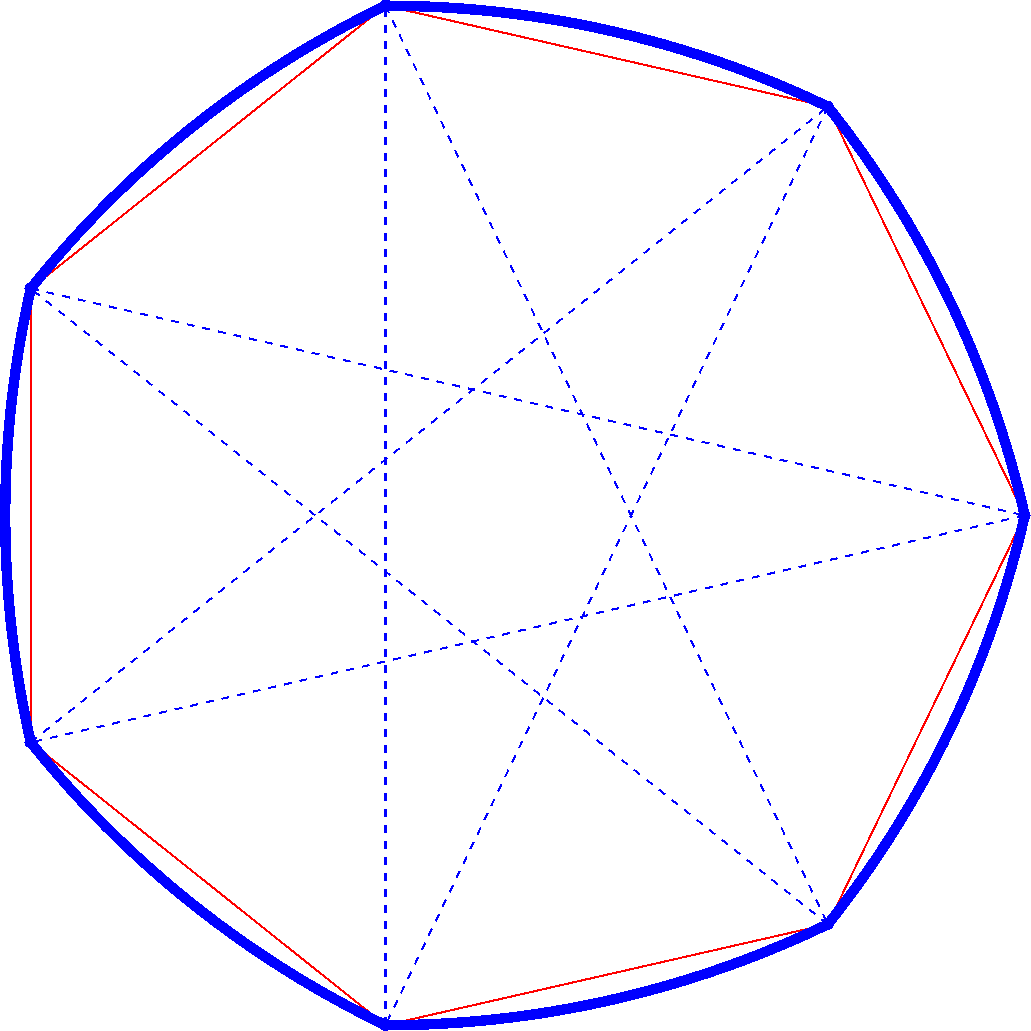}
	\caption{Examples of Reuleaux polygons}
	\label{fig:Rpoly}
\end{figure}

Constant width shapes are convex domains for which any pair of parallel supporting lines are at a fixed distance. In the following, all constant width shapes have unit width. A disk is trivially a shape of constant width, but is not the only one. Intersections of unit disks centered at vertices of a regular $n$-gon with unit diameter generate regular Reuleaux $n$-gons. For $n=3$ the classical Reuleaux triangle is obtained. Examples are shown in Figure \ref{fig:Rpoly}.

The minimization of the area is well posed in the class of constant width shapes and the classical Blaschke-Lebesgue theorem (see \cite{blaschke_orig}, \cite{Lebesgue}) states the following.

\begin{thm}\label{thm:BL}
	The Reuleaux triangle minimizes the area among all shapes of fixed constant width.
\end{thm}

 Many different proofs besides the original ones are known for this result. Geometric arguments are given in the book of Yaglom and Boltyanskii \cite[Chapter 7]{yaglom-boltjanskii}. Chakerian uses circumscribed regular hexagons and mixed volumes in \cite{Chakerian}. Ghandehari uses control theory in \cite{Ghandelhari} and Harrell uses variational techniques in \cite{Harrell}. 
 
 In \cite{Firey1960}, \cite{Sallee70} it is shown that when restricting to the particular class of Reuleaux $n$-gons (the precise definition is given in Section \ref{sec:reuleaux}), the regular one maximizes the area. Kupitz and Martini study in detail isoperimeric inequalities related to Reuleaux $n$-gons in \cite{kupitz-martini}.
 
 The purpose of this note is to find new purely variational arguments, in the same spirit as \cite{Harrell}, for the discrete case of Reuleaux $n$-gons. It turns out that when $n \geq 5$ any vertex of a Reuleaux $n$-gon may be perturbed in a smooth way to generate another Reuleaux $n$-gon. Since Reuleaux polygons of minimal and maximal area exist when considering an upper bound on the number of vertices, optimality conditions of the first and second order can be written for the optimal shapes. 
 
 {\bf Main results.} The optimality conditions show that the only critical points for the area are the regular Reuleaux $n$-gons. An explicit computation shows that among regular Reuleaux $n$-gons, the triangle has the minimal area. A more careful investigation shows that no local minimizers for the area exist among Reuleaux $n$-gons when $n \geq 5$. This observation was also made in \cite{henrot_lucardesi_annulus} using Blaschke perturbations which change only two vertices in a Reuleaux polygon.
 
 Since minimizers and maximizers for the area exist in the class of Reuleaux polygons with at most $n$ vertices (see \cite{kupitz-martini} and Proposition \ref{prop:existenceRP}), an optimizer for the area is either a critical point, i.e. a Regular Reuleaux $n$-gon or the Reuleaux triangle. This alternative comes from the fact that the Reuleaux triangle is the only Reuleaux polygon for which no vertex can be perturbed arbitrarily to generate another Reuleaux $n$-gon without increasing the number of vertices. 
 
 Therefore, the following two conclusions follow from the variational arguments:
 
 (a) The Reuleaux triangle minimizes the area among Reuleaux polygons with at most $n$ vertices. The density of the Reuleaux polygons in the class of constant width shapes implies the Blaschke-Lebesgue theorem.
 
 (b) The regular Reuleaux $n$-gon maximizes the area among Reuleaux polygons with at most $n$ vertices. This result is attributed to Firey \cite{Firey1960} and Sallee \cite{Sallee70}.
 
 \bo{Main motivation: the 3D problem.} The motivation for this variational study for the two dimensional Blaschke-Lebesgue problem is not merely to find different arguments for the two dimensional case. The objective is to find methods which generalize to the three dimensional case. 
 
 The Blaschke-Lebesgue problem in dimension three is still unsolved today. It is conjectured that the 3D analogue of the Reuleaux triangle, the Meissner tetrahedra \cite{Meissner-Schilling}, minimize the volume among shapes of unit constant width. The conjecture is attributed to Bonnesen and Fenchel \cite{Bonnesen-Fenchel}. Historical facts are given in \cite{kawohl-webe} and numerical simulations in \cite{AntunesBogosel22}. Nevertheless, none of the many methods which work in 2D generalize to 3D.
 
 New geometrical objects called Meissner polyhedra, which are the 3D analogue of Reuleaux polygons, were introduced recently in \cite{montejano}. Three dimensional constant width bodies are constructed starting from a finite number of points verifying a combinatorial property. Details and a very nice introductory exposition of these objects is given in \cite{meissner_hynd}. Moreover, the area and volume of these bodies was computed in \cite{hynd-vol-per} and \cite{bogosel_Meissner}. In \cite{bogosel_Meissner} it is shown that Meissner tetrahedra are volume minimizers among Meissner pyramids. 
 
 Nevertheless, like indicated in \cite{bogosel_Meissner}, a proof for the three dimensional case should be found if one understands the variation of the area of Meissner polyhedra with respect to vertex perturbations. A deeper understanding of the two dimensional case is important in this perspective since perturbations of Reuleaux polygons are way simpler than perturbations of Meissner polyhedra, due to combinatorial aspects related to the preservation of diameters between particular vertices. 
 
 \bo{Methods used in the paper.} There are at least two possibilities to study optimality conditions for problems involving Reuleaux polygons. Both are illustrated in this paper and are detailed below.
 
 (a) Investigate the variation of the area for perturbations of Reuleaux polygons. One needs to keep in mind that all such perturbations are non-local, to preserve the constant width constraint. If the sensitivities of the area are found when perturbing one vertex, first and second order classical optimality conditions can be written for an optimal shape.
 
 (b) Consider the Reuleaux polygon as a disk-polygon (introduced in the next section), which does not necessarily have constant width. Then impose the constant width constraint as a series of equality constraints for the distances between the centers of the disks defining the disk-polygon. First and second order optimality conditions can be studied in this framework using the Lagrangian \cite[Chapter 12]{Nocedal_Wright}.
 
 {\bf Structure of the paper.} Basic facts about Reuleaux $n$-gons and disk polygons are recalled in Section \ref{sec:reuleaux}. The computation of the sensitivity of the area with respect to vertex perturbations for disk polygons and Reuleaux polygons can be found in Section \ref{sec:sensitivity}. The Lagrangian perspective is presented in Section \ref{sec:Lagrangian} taking again two points of view: vertex perturbations and center perturbations for disk-polygons.

\section{Reuleaux polygons and disk polygons}
\label{sec:reuleaux}

In this article all shapes of constant width are assumed to have unit width. A particular class of planar constant width shapes are the Reuleaux polygons. The simplest example is the Reuleaux triangle, obtained as the intersection of three unit disks centered at the vertices of an equilateral triangle with edge length one. More precisely let us consider the following definition.

\begin{deff}\label{def:Reuleaux}
	A Reuleaux polygon is a constant width shape given by the intersection of a finite number of unit disks. The minimal number $n\geq 3$ of disks involved in the definition of a Reuleaux polygon also gives the number of \emph{sides} and \emph{vertices}.
\end{deff}

In the following it is assumed that the disks which determine a Reuleaux polygon have pairwise distinct centers. Let us enumerate some basic properties of Reuleaux polygons. 

\begin{prop}\label{prop:Rpoly}
	(a) Without loss of generality it can be assumed that a Reuleaux polygon is an intersection of an odd number $n=2k+1\geq 3$ of unit disks $B(c_i)$, $i=0,...,n-1$. 
	
	(b) Each center of the disks $B(c_i)$ generating a Reuleaux polygon $R$ is an angular point of the boundary of $R$. Denoting by $\theta_i$ the length of the arc $B(c_i) \cap \partial R$ we have $\sum_{i=0}^{n-1} \theta_i=\pi$. Moreover, the lengths of the arcs verify $\theta_i \in [0,\pi/3]$.
\end{prop}

\emph{Proof:} (a) Let $c_0,...,c_{n-1}$ be the distinct vertices of a Reuleaux polygon in counter clockwise order on the boundary of $R$. To each vertex $c_i$ associate the opposite arc $\gamma_i$. The arcs $\gamma_i$ follow the same orientation as the vertices $c_i$ on $\partial R$. Consider $c_k,c_{k+1}$ the endpoints of $\gamma_0$. Then $\gamma_1,...,\gamma_{k}$ are the only arcs of $\partial R$ in the half plane determined by $c_0c_{k+1}$ which does not contain $c_1$. Thus, $\partial R$ is made of the arcs $\arc{c_ic_{i+1}}$, $i=0,...,k$ and $\gamma_1,...,\gamma_k$, giving a total of $n=2k+1$ arcs. 

(b) These properties are classical and geometrical proofs can be found in \cite[Chapter 7]{yaglom-boltjanskii}. \hfill $\square$

 It is well known that Reuleaux polygons can approximate arbitrarily well any planar shape of constant width. For a proof see \cite[Chapter 7]{yaglom-boltjanskii}. This density property makes Reuleaux polygons appropriate for studying optimization problems, since an inequality valid for Reuleaux polygons will extend to all constant width shapes. Some proofs of Theorem \ref{thm:BL} simply show that the area of any Reuleaux polygon is greater than the one of the Reuleaux triangle. See \cite[Chapter 7]{yaglom-boltjanskii} and the recent article \cite{hynd-BL}, where additional details are given.

More precise inequalities can be formulated for areas of Reuleaux polygons. In \cite{Firey1960}, \cite{Sallee70} it is shown that among Reuleaux polygons with a fixed number $n$ of sides the regular Reuleaux polygon, having centers at the vertices of a regular $n$-gon of unit diameter, maximizes the area. This result is discussed in detail in \cite{kupitz-martini}. Let us recall two results.

\begin{prop}\label{prop:existenceRP}
	The area functional admits minimizers and maximizers in the class of Reuleaux polygons with at most $n=2k+1\geq 3$ vertices.
\end{prop}

A proof is given in \cite[Section 2]{kupitz-martini}, but for the sake of completeness a simple short argument is given below.

\emph{Proof:} Assume $n \geq 5$, otherwise there is nothing to prove. A set of constant width has trivial upper and lower bounds for the area (zero and the area of a unit disk). Therefore minimizing and maximizing sequences exist for the area functionals among Reuleaux polygons with at most $n$ sides. The Blaschke selection theorem \cite[Theorem 1.8.7]{schneider} assures the existence of converging sequences in the Hausdorff metric and the area is continuous in this metric. Moreover, since the width (in any direction) is continuous for the Hausdorff convergence, any limit point is also a constant width set.

It remains to show that the family of Reuleaux polygons with at most $n$ vertices/sides is closed under the Hausdorff metric. Consider a minimizing/maximizing sequence $(R_k)$ of Reuleaux polygons with at most $n$ sides. There exists a subsequence, denoted again with $(R_k)$, by simplicity, having exactly $n_0 \leq n$ vertices. 

Up to extracting a diagonal sequence, assume the centers $(c_i^k)$ of $R_k$ converge to $c_i$, for $i=1,...,n_0$. The intersection of disks $B(c_i)$ is again a Reuleaux polygon with at most $n_0 \leq n$ vertices. It is possible that vertices merge at convergence. Nevertheless, it follows that the class of Reuleaux polygons with at most $n$ vertices is closed under the Hausdorff metric. The existence of area minimizers and maximizers follows. \hfill $\square$ 

The following result regarding regular Reuleaux $n$-gons, using an equivalent formula, is also stated in \cite{Firey1960}.

\begin{prop}\label{prop:areaRegR}
	The area of a regular Reuleaux polygon with $n=2k+1\geq 3$ sides is given by
	\begin{equation}\label{eq:areaRegR}
	A_n = \frac{\pi}{2}-n \frac{\sin \frac{\pi}{n}}{2(\cos \frac{\pi}{n}+1)}.
	\end{equation}
	The sequence $A_n$ is strictly increasing.
\end{prop}

\emph{Proof:} Simply observe that a regular Reuleaux $n$-gon is made of a regular $n$-gon with diameter $1$ and $n$ circle sectors having areas $0.5(\frac{\pi}{n}-\sin \frac{\pi}{n})$.  See Figure \ref{fig:Rpoly} for an illustration. This implies that the area of the regular Reuleaux $n$-gon is
\[
A_n = n \frac{\sin(2\pi/n)}{8 \cos^2 (\pi/(2n))}+\frac{1}{2}\left(\pi-n\sin \frac{\pi}{n}\right),
\]
which after a few trigonometric simplifications gives \eqref{eq:areaRegR}.

The monotonicity of $A_n$ is proved in \cite{Firey1960}. Using \eqref{eq:areaRegR} the monotonicity of $A_n$ follows from the fact that 
\[ g: x \mapsto \frac{\sin x}{2x(\cos x+1)}\]
is a strictly decreasing function and $A_n = \frac{\pi}{2}-\pi g(\frac{\pi}{n})$.  \hfill $\square$

Since Reuleaux polygons are intersections of a finite number of unit disks, it will be useful to consider such sets which do not necessarily have constant width.

\begin{deff}\label{def:disk-poly}
	A \emph{disk-polygon} is an intersection of a finite number of disks having unit radius.
\end{deff}

\begin{figure}
	\centering
	\includegraphics[width=0.5\textwidth]{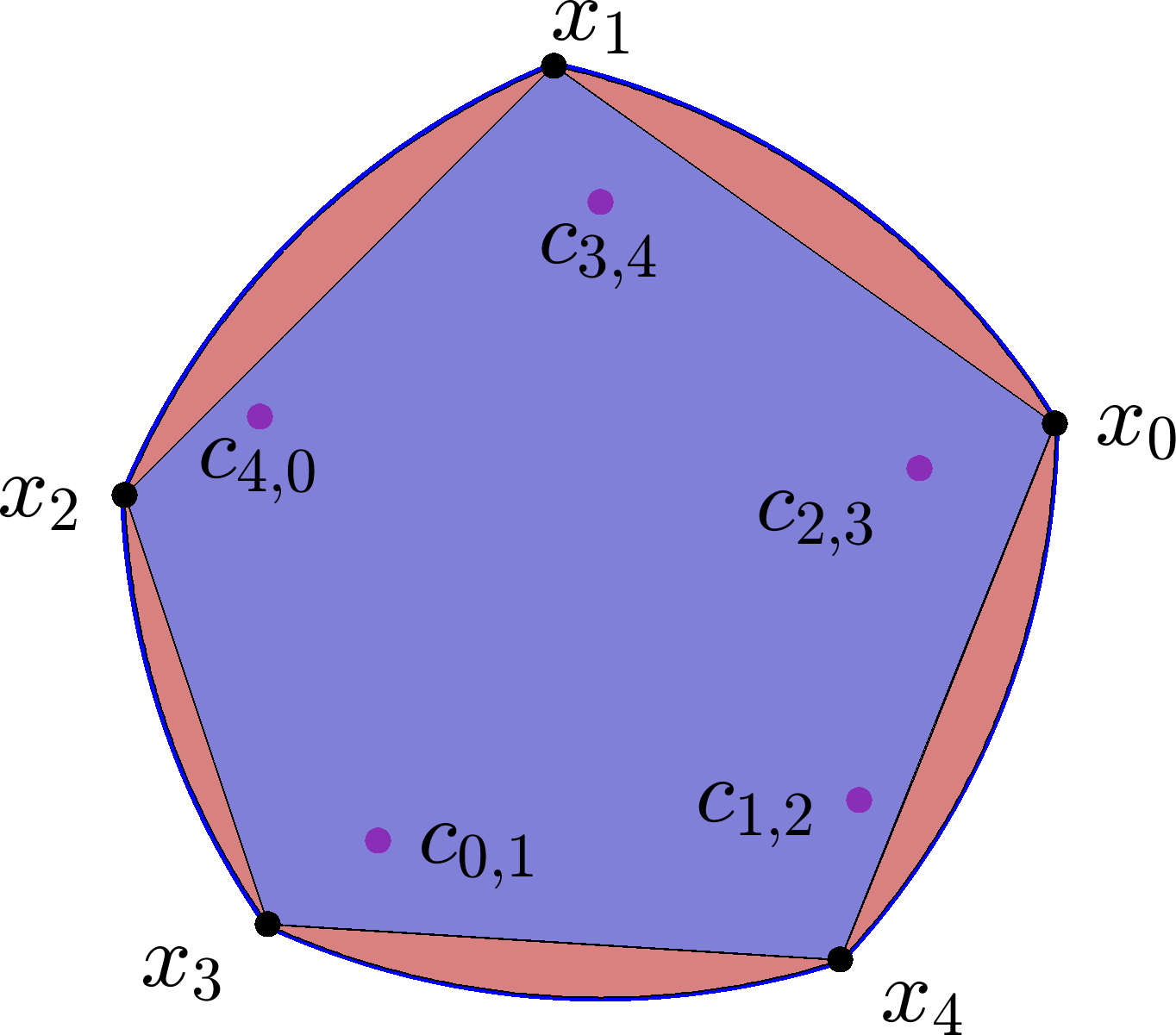}
	\caption{Example of disk polygon with $5$ vertices. The point $c_{i,i+1}$ is the center of the circle containing the arc passing through $x_i$ and $x_{i+1}$ in the boundary of the disk polygon.
	}
	\label{fig:disk-polygon}
\end{figure}

An example is shown in Figure \ref{fig:disk-polygon}. In particular, a Reuleaux polygon is a disk polygon having constant width. Let us recall the following properties:
\begin{itemize}[noitemsep]
	\item If $x_0,...,x_{n-1}$ are the vertices of a disk-polygon $P$ then $P$ is the intersection of all unit disks containing the points $x_0,...,x_{n-1}$.
	\item In general, assuming the vertices are labeled in the anti-clockwise sense, the centers of the disks bounded by the arcs $\arc{x_ix_{i+1}}$ are not necessarily among the vertices of the disk polygon.
	\item A disk polygon can either be defined by specifying the vertices $x_0,...,x_{n-1}$ of a convex $n$-gon in the plane or by specifying the disk centers. In Figure \ref{fig:disk-polygon} the notation $c_{i,i+1}$ is used for the center of the circle containing the arc $\arc{x_ix_{i+1}}$.
\end{itemize}  

\section{Sensitivity of the area with respect to vertex perturbations}
\label{sec:sensitivity}
The study optimization problems related to the areas of Reuleaux polygons is facilitated if one understands how the area is changed when the vertices are perturbed. Many geometric proofs of the Blaschke-Lebesgue theorem show that the area decreases when performing a geometric transformation of the Reuleaux polygon. The aim of this section is to understand the derivatives of the area of Reuleaux polygons and disk polygons in general. In \cite{henrot_lucardesi_annulus} the authors compute the sensitivity of the area for particular perturbations of Reuleaux polygons, called Blaschke perturbations. These perturbations are carefully chosen such that only two vertices move. In the following, the sensitivity of the area will be computed with respect to arbitrary vertex perturbations.

A useful concept when computing the sensitivity of an objective function with respect to boundary variations of a geometric domain is the shape derivative. It is well known that if $V(t)$ is a Lipschitz vector field that perturbs a Lipschitz regular shape $\Omega$ (convex domains, in particular) then the variation of the area is given by the following shape derivative formula 
\begin{equation}\label{eq:area-deriv}
\Area ((I+V(t))(\Omega))'|_{t=0} = \int_{\partial \Omega} V'(0)\cdot \bo n ds,
\end{equation}
where $\bo n$ denotes the outer unit normal.
This result is a consequence of classical shape derivative formulas and may be found, for example, in \cite[Chapter 5]{henrot-pierre-english}. 

The shape derivative formula \eqref{eq:area-deriv} shows that in order to understand the variation of the area, one simply needs to compute the vector field $V$ characterizing the boundary perturbation.

\subsection{Disk-polygons}

Consider a disk-polygon with vertices $x_0,...,x_{n-1}$ (labeled in anti-clockwise order), $n \geq 3$. For each vertex consider a perturbation of the form $x_i(t) = x_i+tv_i$ where $v_0,...,v_{n-1}$ are arbitrary vectors in $\Bbb{R}^2$. Assume each arc $\arc{x_ix_{i+1}}$ is strictly contained in a half unit circle. The shape derivative formula \eqref{eq:area-deriv} is a boundary integral which can be computed separately on each one of the arcs $\arc{x_ix_{i+1}}$. Moreover, the boundary between $x_i(t)$ and $x_{i+1}(t)$ is always an arc of a unit circle, for $t$ small enough, whose center will be denoted by $c_{i,i+1}(t)$. If $V(t)$ is a vector field mapping $\arc{x_ix_{i+1}}$ to $\arc{x_i(t)x_{i+1}(t)}$ for each $t\geq 0$ then the infinitesimal perturbation $V'(0)$ is given by the derivative of the perturbation of the center of the corresponding disk $c_{i,i+1}'(0)$.

\bo{Perturbation of the circle center for a given arc.} The center $c(t)$ verifies for all $t \geq 0$
\[ |x_i+t v_i-c_{i,i+1}(t)| = 1, \ |x_{i+1}+tv_{i+1}-c_{i,i+1}(t)|=1.\]
The mapping $t\mapsto c_{i,i+1}(t)$ is differentiable due to the implicit function theorem, assuming $c_{i,i+1}$ is always on the same side of the segment $x_i(t)x_{i+1}(t)$. Differentiating the above equalities at $t=0$ gives
\[ \begin{array}{rcl}
(x_i-c_{i,i+1}(0))\cdot c_{i,i+1}'(0) &=& (x_i-c_{i,i+1}(0))\cdot v_i\\
(x_{i+1}-c_{i,i+1}(0))\cdot c_{i,i+1}'(0) &=& (x_{i+1}-c_{i,i+1}(0))\cdot v_{i+1}
\end{array}.\]
Since $(x_i-c_{i,i+1}(0)), (x_{i+1}-c_{i,i+1}(0))$ are two linearly independent vectors, the above inequalities completely determine $c_{i,i+1}'(0)$. For simplicity denote 
\begin{equation}\label{eq:notation-w}
 w_i = x_i-c_{i,i+1}(0), w_{i+1}=x_{i+1}-c_{i,i+1}(0)
 \end{equation}
and decompose $c_{i,i+1}'(0) = \alpha_i w_i+\alpha_{i+1} w_{i+1}$. Also denote by $\theta_{i,i+1}\in (0,\pi/2)$ the arclength of $\arc{x_ix_{i+1}}$. It follows that 
\begin{equation}\label{eq:alpha}
 \begin{pmatrix}
1 & \cos\theta_{i,i+1} \\
\cos\theta_{i,i+1} & 1
\end{pmatrix} \begin{pmatrix}
\alpha_i\\ \alpha_{i+1}
\end{pmatrix} = \begin{pmatrix}
v_i \cdot w_i\\
v_{i+1}\cdot w_{i+1}
\end{pmatrix}.
\end{equation}
The following result helps compute the contribution of the arc $\arc{x_ix_{i+1}}$ to the shape derivative of the area of a disk-polygon.

\begin{prop}\label{prop:one-arc}
	a) Let $\gamma$ be an arc of a unit circle of arclength $\theta$. Let $\bo b$ be the unit vector bisecting the corresponding circle sector pointing from the center towards the midpoint of the arc. Then for $v \in \Bbb{R}^2$ an arbitrary fixed vector we have
	\[ \int_{\gamma}  v\cdot \bo n = 2\sin \frac{\theta}{2}  v\cdot \bo b.\]
	
	b) If $v_i,v_{i+1}$ are perturbation vectors for vertices $x_i,x_{i+1}$ of a disk polygon and $V'(0)$ is the corresponding infinitesimal perturbation then
	\[ \int_{\arc{x_ix_{i+1}}} V'(0)\cdot \bo n = \tan \frac{\theta_{i,i+1}}{2}[(x_i-c_{i,i+1})\cdot v_i+(x_{i+1}-c_{i,i+1})\cdot v_{i+1}].\]
\end{prop}

\emph{Proof:} a) First, observe that if $ t$ is orthogonal to $\bo b$ then $ t\cdot \bo  n$ is an odd function on $\gamma$, whose integral is zero. Therefore, since $\bo b \cdot \bo n$ is symmetric on $\gamma$, we have
\[ \int_\gamma  v\cdot \bo n = (v \cdot \bo b)\int_\gamma \bo b\cdot \bo  n=( v \cdot \bo b) 2\int_0^{\theta/2} \cos t dt=2\sin \frac{\theta}{2}  v\cdot \bo  b.\]

b) The previous observations and result of point a) show that
\[ \int_{\arc{x_ix_{i+1}}} V'(0)\cdot  \bo  n =  \int_{\arc{x_ix_{i+1}}} c_{i,i+1}'(0)\cdot \bo n = 2\sin \frac{\theta_{i,i+1}}{2} c'_{i,i+1}(0)\cdot \bo b_{i},\]
where $\bo b_i$ is the unit bisector of the angle $\angle x_ic_{i,i+1}x_{i+1}$.

Using the notations \eqref{eq:notation-w} gives
\[\int_{\arc{x_ix_{i+1}}} V'(0)\cdot \bo n = 2\sin \frac{\theta_{i,i+1}}{2}( \alpha_i w_i\cdot \bo b_{i}+\alpha_{i+1} w_{i+1}\cdot \bo b_{i+1})= 2\sin \frac{\theta_{i,i+1}}{2}\cos \frac{\theta_{i,i+1}}{2}(\alpha_i+\alpha_{i+1}).\]
Finding $\alpha_i,\alpha_{i+1}$ solutions of \eqref{eq:alpha} gives
\[\int_{\arc{x_ix_{i+1}}} V'(0)\cdot \bo n = \tan \frac{\theta_{i,i+1}}{2}(w_i \cdot v_i+w_{i+1}\cdot v_{i+1}),\]
as requested.
\hfill $\square$ 

The shape derivative of the area of a disk polygon $P$ is obtained by summing the contributions in Proposition \ref{prop:one-arc} for each one of the arcs in the boundary of $P$. The previous arguments remain valid if $x_i, x_{i+1}$ are perturbed using more general maps $t \mapsto x_i(t)$ such that $x_i'(0)=v_i$, $x_{i+1}'(0)=v_{i+1}$.

\subsection{Reuleaux polygons}
The shape derivative of the area of a Reuleaux polygon with respect to vertex perturbations will be found using similar ideas to disk-polygons. The key difference is that, while for disk-polygons perturbations $v_i$ for vertices $x_i$, $i=0,...,n-1$ may be arbitrary, for Reuleaux polygons the vertex perturbations must verify certain constraints which preserve the constant width. 

Consider a Reuleaux polygon $R$ with $n=2k+1\geq 5$ sides and vertices $x_0,...,x_{n-1}$. The length of the arc opposite to $x_i$ is denoted by $\theta_i$. As stated in Definition \ref{def:Reuleaux} all centers of the disks generating $R$ are pairwise distinct. 

In the following, it is shown that every vertex of $R$ may be perturbed, generating another Reuleaux $n$-gon. It is well known that the angles/lengths of the arcs in $\partial R$ verify  $\theta_i \in [0,\pi/3]$ \cite[Exercises 7-9]{yaglom-boltjanskii}. The result follows at once noticing that any disk sector of opening $\theta_i$ must have diameter $1$. Moreover, if $\theta_i=\pi/3$ for some $i$ then $R$ is a Reuleaux triangle. Thus, if $R$ is a Reuleaux $n$-gon with $n \geq 5$, all angles verify $0<\theta_i < \pi/3$, $i=0,...,n-1$.

 Consider the vertex $x_0$ of $R$ and an arbitrary unit vector $ v \in \Bbb{R}^2$. Assume the vertices $x_0,...,x_{n-1}$ are oriented  counter-clockwise on the boundary of $R$. Then the opposite edge to $x_0$ is $\arc{x_kx_{k+1}}$. Moreover, $x_k \in \partial B(x_{n-1})$ and $x_{k+1}\in \partial B(x_1)$. 

\begin{figure}
	\centering 
	\includegraphics[height=0.4\textwidth]{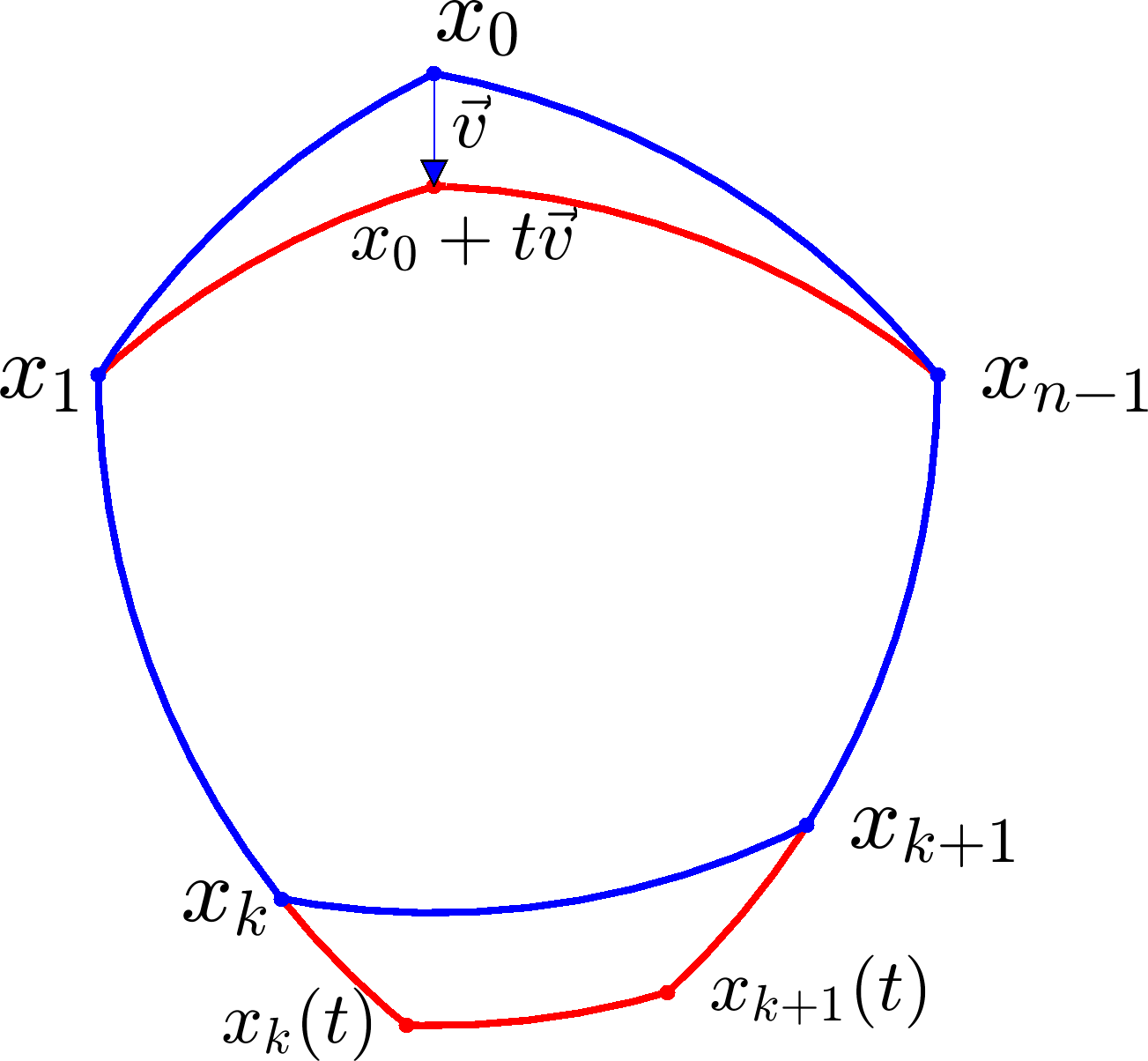}
	\includegraphics[height=0.4\textwidth]{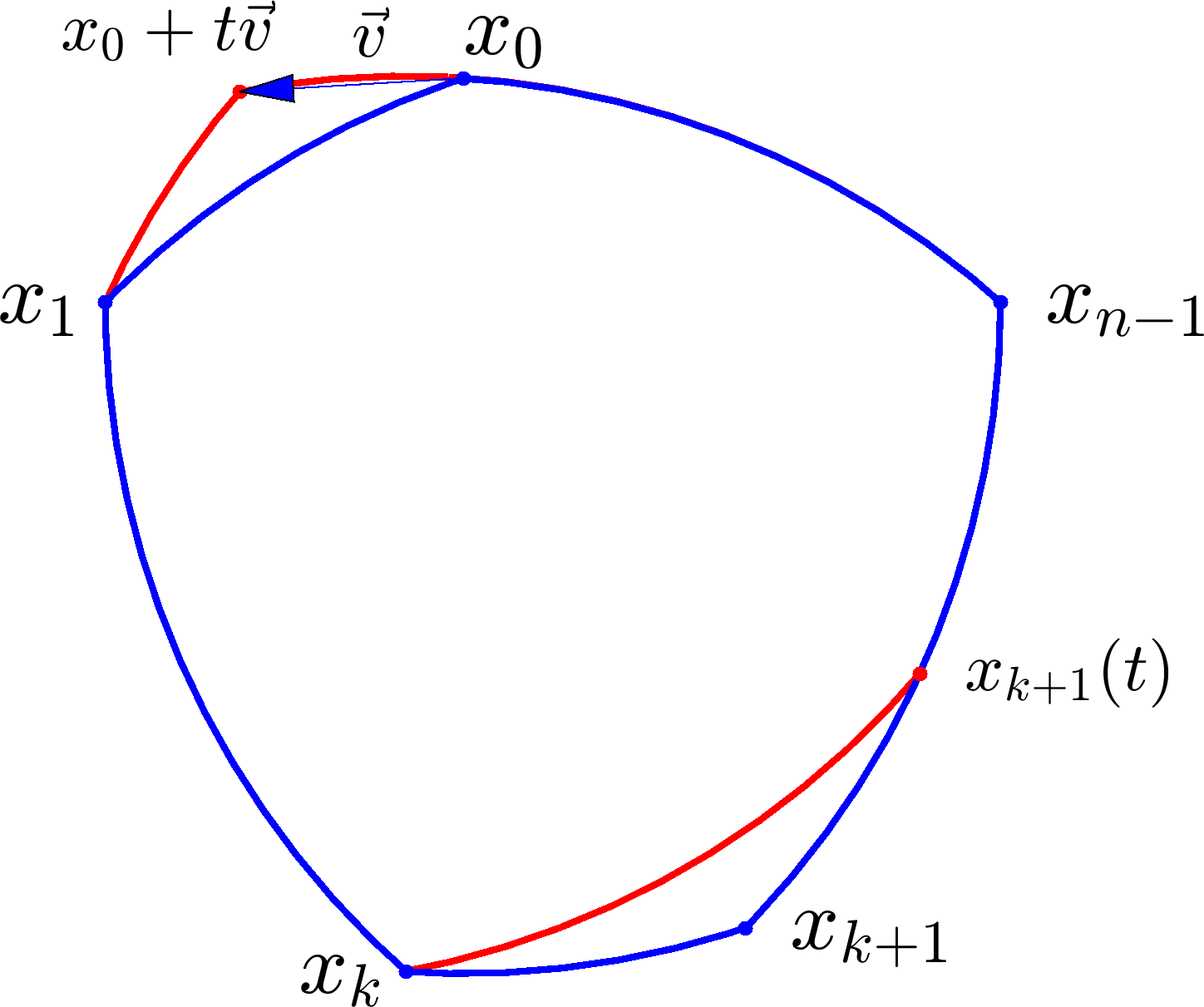}
	\caption{General perturbation of a vertex of a Reuleaux polygon generating a new Reuleaux $n$-gon. The left picture shows a general perturbation. The right picture shows a Blaschke perturbation: the vertex $x_0$ is perturbed tangentially, leading to a simpler boundary variation: only two vertices move.}
	\label{fig:perturbedR}
\end{figure}

For $t$ small enough, consider 
\begin{equation}\label{eq:perturbed-centers}
x_k(t) = \partial B(x_0+t v)\cap \partial B(x_{n-1}) \text{ and }  x_{k+1}(t) = \partial B(x_0+t v) \cap \partial B(x_1).
\end{equation} 
Thus a new perturbed Reuleaux $n$-gon is obtained, denoted by 
\[ R(t) = B(x_0+t v)\cap B(x_1)\cap ... \cap B(x_k(t))\cap B(x_{k+1}(t)) \cap ... \cap B(x_{n-1}).\]
See Figure \ref{fig:perturbedR} for an illustration. Indeed, since $\theta_i<\pi/3$ for all $i=1,...,n$ and the transformation $t\mapsto x_0+t v$ is continuous, the same inequality $\theta_i(t)<\pi/3$ holds for $t>0$ small enough, where $\theta_i(t)$ are the edge lengths of $R(t)$. Thus $R(t)$ is still a Reuleaux $n$-gon. Notice that this type of perturbation cannot be made for a Reuleaux triangle. 

In the following, the sensitivity of the area of $R(t)$ with respect to $t$ is computed. The fact that $R(t)$ is a Reuleaux polygon follows from the construction and from \eqref{eq:perturbed-centers}. Depending on the orientation of $ v$ the arcs $\arc{x_0x_1}$, $\arc{x_0x_{n-1}}$, $\arc{x_{k-1}x_k}$, $\arc{x_kx_{k+1}}$, $\arc{x_{k+1}x_{k+2}}$ may be modified. Initially all these arcs have lengths in $(0,\pi/3)$ for $t$ small enough. Since the perturbation is continuous, the perturbed arcs will verify the same inequalities.

 When moving vertices $x_0, x_k, x_{k+1}$ in the definition of $R(t)$ normal displacements of the boundary occur only on the arcs $\arc{x_0x_1}, \arc{x_0x_{n-1}}, \arc{x_kx_{k+1}}$. The following result follows by applying Proposition \ref{prop:one-arc}.

\begin{thm}\label{thm:sensitivity}
	Suppose $ v$ is a unit vector and denote $\bo b$ the unit vector aligned with the angle bisector of $\angle x_kx_0x_{k+1}$.
	
	(a) The directional derivative of the area of $R(t)$ when $x_0$ is perturbed in the direction $v$ is given by
	\begin{equation}\label{eq:sensitivity-one-vertex}
	|R(t)|'_{t=0} =2\sin \frac{\theta_0}{2} ( v\cdot \bo b) + \tan \frac{\theta_{k}}{2} v\cdot (x_0-x_k)+\tan \frac{\theta_{k+1}}{2} v\cdot (x_0-x_{k+1}).
	\end{equation}
    If $|R(t)|_{t=0}'=0$ for all directions $ v$ then $\theta_0 = \theta_k = \theta_{k+1}$. In particular, $\theta_0,\theta_k,\theta_{k+1} \in [0,\pi/5]\cup \{\pi/3\}$. 
	
	(b) If $\theta_0=\theta_k=\theta_{k+1} \leq \pi/5$ then $R$ is not a local minimizer for the area with respect to perturbations of $x_0$.
\end{thm}

\emph{Proof:} (a) The center of $\arc{x_0x_1}$ is $x_{k+1}$ and its length is $\theta_{k+1}$. Similarly, the center of $\arc{x_{n-1}x_0}$ is $x_k$ and its length is $\theta_k$.

In view of Proposition \ref{prop:one-arc} if $x_0$ is perturbed by $x_0(t) = x_0+tv$ and vertices $x_1,x_{n-1}$ are not moved then the contribution to arcs $\arc{x_{n-1}x_0},\arc{x_0x_1}$ in the shape derivative formula \eqref{eq:area-deriv} is given by 
\[ \tan\frac{\theta_k}{2} (x_0-x_k)\cdot v +\tan \frac{\theta_{k+1}}{2}(x_0-x_{k+1})\cdot v.\]

Perturbing $x_0$ in direction $v$ induces the normal infinitesimal perturbation $v\cdot n$ on $\arc{x_kx_{k+1}}$, which is an arc of circle centered at $x_0$. Thus, the contribution of this arc to the shape derivative is
\[ \int_{\arc{x_kx_{k+1}}} v\cdot\bo n = 2\sin \frac{\theta_0}{2} (v\cdot \bo  b),\]
according to Proposition \ref{prop:one-arc} a). Combining the contributions given by the arcs $\arc{x_{n-1}x_0},\arc{x_0x_1},\arc{x_kx_{k+1}}$ gives \eqref{eq:sensitivity-one-vertex}.


Suppose that $R$ is a critical point for the area when considering perturbations of $x_0$. Then \eqref{eq:sensitivity-one-vertex} implies
\[ 2\sin \frac{\theta_0}{2} \bo  b = -\tan \frac{\theta_{k+1}}{2} (x_0-x_{k+1})-\tan \frac{\theta_k}{2} (x_0-x_k).\]
Since $\bo b$ is the bisector of $\angle x_kx_0x_{k+1}$ and $|x_0-x_k|=|x_0-x_{k+1}|=1$ it follows that
\begin{equation}\label{eq:bisector}
\bo b = \frac{1}{2\cos \frac{\theta_0}{2}}(x_k-x_0+x_{k+1}-x_0).
\end{equation}
Therefore, the following vectorial equality holds
\[\tan \frac{\theta_{0}}{2} \overrightarrow{x_0x_{k+1}}+\tan \frac{\theta_0}{2} \overrightarrow{x_0x_k} = \tan \frac{\theta_{k+1}}{2} \overrightarrow{x_0x_{k+1}}+\tan \frac{\theta_k}{2} \overrightarrow{x_0x_k},\]
which implies $\theta_0=\theta_k=\theta_{k+1}$.

To observe the bounds on $\theta = \theta_0=\theta_k=\theta_{k+1}$ consider the length of $x_1x_{n-1}$ in terms of $\theta$. A simple computation using the law of cosines in $\Delta x_0x_1x_{n-1}$ gives
\[ x_1x_{n-1} = -8\sin^3\frac{\theta}{2}+4\sin \frac{\theta}{2}.\]
The diameter constraint implies $x_1x_{n-1}\leq 1$, thus we are led to the inequality 
\[-8\sin^3\frac{\theta}{2}+4\sin \frac{\theta}{2}-1 \leq 0.\]
Factoring we obtain
\[ \left(1-2\sin \frac{\theta}{2}\right)\left(4\sin^2 \frac{\theta}{2}+2\sin \frac{\theta}{2}-1\right) \leq 0.\]
The inequality is satisfied when $\theta = \pi/3$ or when $\sin \frac{\theta}{2} \leq \frac{\sqrt{5}-1}{4} = \sin \frac{\pi}{10}$. Therefore at any critical configuration the angles satisfy $\theta \in [0,\pi/5]\cup \{\pi/3\}$.

(b) The formula for the first derivative of $|R|$ when perturbing $x_0$ in direction $ v$ depends on $\theta_0,\theta_k, \theta_{k+1}$ and the vector $v$. These angles have directional derivatives, as shown below. Consider a unit direction $ w$ (potentially different than $ v$). 

Denote $\ell_{k+1}(t) = |x_0+t w-x_1|$. Then it is straightforward to see that
\[ \ell_{k+1}'(0) = -\cos(\angle  w, \overrightarrow{x_0x_1}).\]
Since $\theta_{k+1}(t) = 2\arcsin(\ell_{k+1}(t)/2)$ we find
\[ \theta_{k+1}'(0) = -\frac{\cos(\angle  w,\overrightarrow{x_0x_1})}{\cos \frac{\theta_{k+1}}{2}}.\]
Similarly 
\[ \theta_{k}'(0) = -\frac{\cos(\angle  w,\overrightarrow{x_0x_{n-1}})}{\cos \frac{\theta_{k}}{2}}.\]

Assume, for simplicity that $ w$ points inside the angle $\angle x_kx_0x_{k+1}$. Denoting $\ell_0(t)$ the length of the segment ${x_k(t)x_{k+1}(t)}$ we have
\begin{align*} \ell_0'(0) &= - |x_k'(0)|\cos(\angle( x_k'(0),\overrightarrow{x_{k}x_{k+1}}))-|x_{k+1}'(0)|\cos(\angle (x_{k+1}'(0),\overrightarrow{x_{k+1}x_{k}}))\\
&= -\frac{ \cos (\angle (\overrightarrow{x_1x_{k}}, w))}{\sin \theta_{k}}\cos\left(\theta_{k+1}+\frac{\theta_1}{2}\right)+\frac{ \cos (\angle (\overrightarrow{x_1x_{k+1}}, w))}{\sin \theta_{k+1}}\cos\left(\theta_{k}+\frac{\theta_1}{2}\right)
\end{align*}
Since $\theta_0(t) = 2\arcsin \frac{\ell_0(0)}{2}$ it follows that
\[ \theta_0'(0) = \frac{1}{\sqrt{1-\ell_0(0)^2/4}}\ell_0'(0)= \frac{1}{\cos \frac{\theta_0}{2}}\ell_0'(0).\]

Suppose that $R$ is a critical point with respect to perturbations of $x_0$. Therefore $\theta_0=\theta_k=\theta_{k+1}$. Then for $ v$ bisector of the angle $\angle x_kx_0x_{k+1}$, using the symmetry, we have 
\[ R'(t) = 2\sin \frac{\theta_0(t)}{2}-2\cos \frac{\theta_0(t)}{2}\tan \frac{\theta_k(t)}{2}.\]
Differentiating again with respect to $t$ at $0$ gives
\[ R''(0) = \left(\cos \frac{\theta_0}{2}+\sin \frac{\theta_0}{2}\tan \frac{\theta_k}{2}\right)\theta_0'(0)-\cos \frac{\theta_0}{2} \frac{1}{\cos^2(\theta_k/2)}\theta_k'(0).\]
Using $\theta_0=\theta_k$ and the formulas found above for $\theta_0'(0), \theta_k'(0)$ we find
\[ R''(0) =2 \frac{\sin \frac{\theta}{2} \sin \theta - \cos \frac{3\theta}{2}}{\cos \frac{\theta}{2}\sin \theta}.\]
Assuming $\theta \in [0,\pi/3]$, the inequality $R''(0)\leq 0$ is equivalent to $\tan\frac{\theta}{2}\tan \theta \leq \frac{1}{2}$. This is obviously verified whenever $\theta \in [0,\pi/5]$, i.e. for any critical configuration excepting the Reuleaux triangle. Since $R''(0)<0$, the area is strictly concave for perturbations along the bisector of $\angle x_kx_0x_{k+1}$. Therefore a critical Reuleaux polygon $R$ with more than $5$ sides is not a local minimizer for the area. \hfill $\square$


The previous result implies that the only critical points for the area among Reuleaux polygons are the regular ones. Combining this with Proposition \ref{prop:areaRegR} gives the following result.

\begin{cor}\label{cor:BL-1}
	Among Reuleaux polygons with at most $n=2k+1\geq 3$ vertices:
	
	(a) the regular Reuleaux $n$-gon maximizes the area.
	
	(b) the Reuleaux triangle minimizes the area.
	
	(c) a Reuleaux polygon which is not a Reuleaux triangle is not a local minimizer for the area. Moreover, any vertex can be perturbed, leaving the other unconstrained vertices fixed, such that the area is strictly decreased.
\end{cor}

This result is also underlined in \cite[Corollary 7]{henrot_lucardesi_annulus} using a family of particular perturbations moving only two vertices of a Reuleaux polygon, called Blaschke perturbations (see Figure \ref{fig:perturbedR}).

\emph{Proof:} (a) Consider the Reuleaux polygon  $R$ with at most $n$ sides which maximizes the area. The existence of a maximizer follows from Proposition \ref{prop:existenceRP}. If $n=3$ there is nothing to prove, since there exists only one Reuleaux triangle. 

If $n\geq 5$ then any one of the vertices can be perturbed as in Theorem \ref{thm:sensitivity}. Since $R$ maximizes the area, it is a critical point when perturbing each one of the vertices, which by Theorem \ref{thm:sensitivity} shows that $R$ is regular. It is not known \emph{a priori} if $R$ has exactly $n$ sides. Nevertheless, since the area of regular Reuleaux $n$-gons is increasing with respect to $n$ (Proposition \ref{prop:areaRegR}), a maximizer must have exactly $n$ sides. 

(b) A minimizer for the area $R$ exists in the class of Reuleaux polygons with at most $n$ sides according to Proposition \ref{prop:existenceRP}. If $R$ has $n \geq 5$ sides then Theorem \ref{thm:sensitivity} and the criticality of $R$ implies that $R$ is regular. On the other hand, Proposition \ref{prop:areaRegR} shows that the Reuleaux triangle has minimal area among regular Reuleaux polygons. Therefore $R$ is not optimal if $n \geq 5$. It follows that the Reuleaux triangle minimizes the area. 

(c) Suppose $R$ is a Reuleaux polygon with $n\geq 5$ vertices which is a local minimzier for the area. Then $R$ is a critical point, and Theorem \ref{thm:sensitivity} implies that $R$ is not a local minimizer, contradicting the assumption.

Consider variations for one vertex  of $R$, denoted with $x_0$ after an eventual relabeling. The constant width constraint implies that vertices $x_k$, $x_{k+1}$ will also be perturbed when moving $x_0$. If $R$ is not a critical point for the area when perturbing $x_0$, then there exists a perturbation of $x_0$ which strictly decreases the area. If $R$ is a critical point for perturbations of $x_0$ then the adjacent and opposite arcs to $x_0$ have equal lengths: $\theta_0=\theta_k=\theta_{k+1}$. Theorem \ref{thm:sensitivity} shows that the area is concave when perturbing $x_0$ in the direction of the bisector of $\angle x_kx_0x_{k+1}$. Therefore the area of $R$ is decreased when perturbing $x_0$ along this direction. 
 \hfill $\square$
 
In \cite{henrot_lucardesi_annulus} the authors define Blaschke perturbations as perturbations of the vertex $x_0$ tangential to one of the neighboring arcs $\arc{x_0x_{n-1}}$, $\arc{x_0x_1}$ (see Figure \ref{fig:perturbedR}). In \cite[Proposition 2.5]{henrot_lucardesi_annulus} the authors give an expression of the derivative of the area for these Blaschke perturbations. Let us see the expression of the derivative for Blaschke perturbations given by Theorem \ref{thm:sensitivity}.
 
 \begin{cor}\label{cor:BL-2} (Sensibility of the area for Blaschke perturbations)
 	
 	(a) If $x_0$ moves on the arc $\arc{x_0x_{n-1}}$ such that the length of this arc increases then the derivative of the area equals
 	\[ R'(0) = \sin \theta_0 \left(\tan \frac{\theta_{k+1}}{2}-\tan \frac{\theta_0}{2} \right).\]
 	
 	(b) The area is concave along a Blaschke perturbation described above when $\theta_0 \geq \theta_{k+1}$. Moreover, if $\theta_{k+1}=\theta_0$ then the configuration gives a local maximum for the area.
 \end{cor}
 
 \emph{Proof:} (a) The value of the derivative of the area is obtained From Theorem \ref{thm:sensitivity} noting that if $\vec v$ is tangent to $\arc{x_0x_{n-1}}$ then $ v \cdot \bo b = -\sin \frac{\theta_0}{2}$, $ v\cdot \overrightarrow{x_0x_{k+1}} = -\sin \theta_0$, $\vec v \cdot \overrightarrow{x_0x_{k}}=0$.
 
 (b) By definition, Blaschke perturbations which increase the arc $\arc{x_0x_{n-1}}$ do the following: $\theta_0$ increases, $\theta_k$ increases, $\theta_{k+1}$ decreases (see Figure \ref{fig:perturbedR}, right). If $\theta_0\geq \theta_{k+1}$ then $\sin \theta_0$ increases and $\tan \frac{\theta_0}{2}-\tan \frac{\theta_{k+1}}{2}$ is increasing. Thus $R'(0)<0$ and $|R'(0)|$ is increasing, showing that $R'(0)$ is decreasing. Thus, the area is concave along a Blaschke perturbation. 
 
 If $\theta_0 = \theta_{k+1}$ along a Blaschke perturbation, since $\theta_0$ increases and $\theta_{k+1}$ decreases it follows that $R'(0)$ is positive before this instance and $R'(0)$ is negative afterwards. This indicates that when $\theta_0 = \theta_{k+1}$ a local maximum occurs. 
 \hfill $\square$ 
 
 \begin{rem}
 	Classical proofs of the Blaschke-Lebesgue theorem relying on Blaschke perturbations are, in fact, a consequence of Corollary \ref{cor:BL-2}. A carefully chosen Blaschke perturbation will decrease the area. Performing the Blaschke perturbation until two vertices merge will decrease the number of sides of a Reuleaux polygon. The procedure is repeated, decreasing the area, until the Reuleaux triangle is reached.
 \end{rem}
 
 Theorem \ref{thm:sensitivity} gives an explicit formula for the directional derivative which is written as a scalar product with the vertex perturbation. Therefore, the gradient of the area with respect to the coordinates of a vertex may be found.
 
 \begin{cor}\label{cor:gradient}
 	The gradient of the area of a Reuleaux polygon having $n=2k+1$ sides with respect to the perturbation of the vertex $x_0$ is given by
 	\[ \frac{\partial |R|}{\partial x_0} = \left(\tan\frac{ \theta_k}{2}-\tan \frac{\theta_0}{2}\right) (x_0-x_k) + \left(\tan \frac{\theta_{k+1}}{2}-\tan \frac{\theta_0}{2}\right)(x_0-x_{k+1}).\]
 \end{cor}

 The result follows from Theorem \ref{thm:sensitivity} and the expression of the bisector vector \eqref{eq:bisector}. This formula shows again that the only critical Reuleaux polygons for the area are the regular ones. Observe that the expression of the gradient with respect to $x_0$ found above depends on the lengths of the arcs $\arc{x_{n-1}x_0},\arc{x_0x_1},\arc{x_kx_{k+1}}$ and the vectors $x_0-x_k, x_0-x_{k+1}$, the diameters reaching the vertex $x_0$.
 
 \section{Second perspective: Lagrangian formulation}
 \label{sec:Lagrangian} 
 
 In the previous sections the vertices of the Reuleaux polygon were constrained to preserve the constant width property. In particular, if a vertex is perturbed, the vetrices on the opposite arc move accordingly on arcs of circles (see Figure \ref{fig:perturbedR}). This restriction is limiting and poses difficulties when computing the first and second derivatives of the area, since the constant width constraint must be preserved by the perturbation field. In this section a Lagrangian framework is proposed, allowing more general perturbations for vertices or centers, while adding Lagrange multipliers for the constant width constraint. The advantage is that Reuleaux polygons are replaced with disk-polygons described by $n=2k+1$ free points in $\Bbb{R}^2$ and the constant width constraint is imposed with a series of equality and inequality constraints regarding the centers of these disks (see Definition \ref{def:disk-poly}). First and second order optimality conditions can still be written for the minimization of the area at the price of using Lagrange multipliers.
 
 \subsection{Vertices as variables}
 
 It is possible to write the Blaschke-Lebesgue problem exclusively in terms of vertices of the disk-polygon as generic points in $\Bbb{R}^2$, verifying a number of constraints. Like before, take $n=2k+1$, $k \geq 1$ an odd positive integer. Let $x_0,...,x_{n-1}$ be the vertices of a convex polygon in $\Bbb{R}^2$. Then denote by $D(x_0,...,x_{n-1})$ the area of the disk-polygon determined by points $x_i$, $i=0,...,n-1$, i.e., the intersection of all unit disks containing the points $x_i$. The area $D(x_0,...,x_{n-1})$ is composed of the area of the polygon $x_0,...,x_{n-1}$ and the area of unit disk segments corresponding to chords of lengths $|x_i-x_{i+1}|$. See Figure \ref{fig:disk-polygon} for an illustration.
 
 The subset of $\Bbb{R}^2$ corresponding to vertices of Reuleaux $n$-gons is defined by
 \begin{equation}\label{eq:constr}
X_n = \{(x_i)_{i=0}^{n-1} :  |x_i-x_{i+k}|=1 \text{ and } |x_i-x_j|\leq 1, \forall 1\leq i,j\leq n.\}
 \end{equation}
 Recall that $k$ is given by $2k+1=n$.
 
  The Blaschke Lebesgue problem becomes
  \begin{equation}\label{eq:area_R_coords}
  \min_{(x_0,...,x_{n-1})\in X_n} |D(x_0,...,x_{n-1})|.
  \end{equation}
  where $|\cdot |$ denotes the area function. Since the area is translation invariant, fixing one point in $X_n$ and leaving the other points variable, verifying the constraints, renders $X_n$ compact. Therefore problem \eqref{eq:area_R_coords} has solutions, providing an alternative proof of existence for the Blachke-Lebesgue problem for Reuleaux polygons having at most $n$ vertices than the one shown in Proposition \ref{prop:existenceRP}.
  
  From now on, assume that $\bo x^* = (x_0,...,x_{n-1})$ is a solution of \eqref{eq:area_R_coords}. Therefore $\bo x^*$ contains the vertices of a Reuleaux polygon; notations from previous sections are used in the following. The theory of constrained optimization problems, recalled for example in \cite[Chapter 12]{Nocedal_Wright}, allows to write optimality conditions with the aid of Lagrange multipliers. 
  
  \bo{First order optimality conditions.} Since the constraint set $X_n$ is defined by a series of inequality conditions, among which a precise subset are verified with equality, we have the following:
  
  \begin{itemize}[noitemsep]
  	\item The objective function and the constraints are defined on $\Bbb{R}^{2n}$ and take real values. These functions are of class $C^1$.
  	\item Constraints $\mathcal C_i(\bo x) =: |x_i-x_{i+k}|-1$ are active, in the sense that $\mathcal C_i(\bo x^*)=0$ at the optimum.
  	\item The constraint functions are of class $C^1$ on $X_n$ since the norm is evaluated on non-zero vectors. Gradients of the active constraints at the optimum are given by
  	\[ \partial_{x_j} \mathcal C_i(\bo x^*) = \begin{cases}
  	\frac{x_i-x_{i+k}}{|x_i-x_{i+k}|} & j=i\\
  	\frac{x_{i+k}-x_{i}}{|x_i-x_{i+k}|} & j=i+k\\
  	0 & \text{ otherwise }
  	\end{cases}.\]
  	The full gradients $\nabla \mathcal C_i(\bo x^*)=(\partial_{x_j} \mathcal C_i(\bo x^*))_{j=0}^{n-1}$ are linearly independent.
  \end{itemize}
Therefore, there exist Lagrange multipliers $\lambda_i$, $i=0,...,n-1$ such that for all $j=0,...,n-1$ we have
\begin{equation}\label{eq:lagrange-first} \partial_{x_j} |D(x_0,...,x_{n-1})|+\sum_{i=0}^{n-1} \lambda_i \partial_{x_j} \mathcal C_i(\bo x^*) = 0.
\end{equation}

Perturbing the vertex $x_j$ in the disk polygon $D(x_0,...,x_{n-1})$ only changes the adjacent arcs $\arc{x_{j-1}x_j}$, $\arc{x_jx_{j+1}}$ having lengths $\theta_{j+k}, \theta_{j-k}$, respectively and being contained into unit circles of centers $x_{j+k}$, $x_{j-k}$, respectively. Therefore, Proposition \ref{prop:one-arc} implies that
\[ \partial_{x_j} |D(x_0,...,x_{n-1})|= \tan \frac{\theta_{j+k}}{2}(x_j-x_{j+k})+\tan \frac{\theta_{j-k}}{2}(x_j-x_{j-k}). \]
The optimality condition becomes
\[  \tan \frac{\theta_{j+k}}{2}(x_j-x_{j+k})+\tan \frac{\theta_{j-k}}{2}(x_j-x_{j-k})+\lambda_j (x_j-x_{j+k}) +\lambda_{j-k}(x_j-x_{j-k})=0,\]
for every $j$, $0 \leq j \leq n-1$. Since $x_j-x_{j+k}$ and $x_j-x_{j-k}$ are linearly independent for every $j$, direct identification gives
\[ \lambda_j = -\tan \frac{\theta_{j+k}}{2},\ \ \lambda_{j-k} = -\tan \frac{\theta_{j-k}}{2}.\]
Observe the index pairings for the Lagrange multipliers and the corresponding angle function: $j\mapsto j+k$ and $j-k\mapsto j-k$.

The previous relation not only gives the explicit expression of the Lagrange multipliers in terms of the arclengths $\theta_j$, $0 \leq j \leq n-1$, but also shows that all multipliers are equal and therefore the optimizer $\bo x^*$ is a regular Reuleaux $n$-gon. This result coincides with the one shown in Corollary \ref{cor:BL-1}, but the proof is more straightforward, since computing area sensitivity for disk-polygons is less restrictive than computing area sensitivity for Reuleaux polygons.

\bo{Second order optimality conditions.} Assume that $\bo x^*$ corresponds to a critical point with at least $5$ vertices. Therefore $(x_0,...,x_{n-1})$ are the vertices of a regular Reuleaux $n$-gon and $\theta_j = \pi/n$, $0 \leq j \leq n-1$. Previous computations show that the Lagrange multipliers appearing in the first order optimality conditions are given by $\lambda_j = -\tan \frac{\pi}{2n}$, $0 \leq j \leq n-1$.

Like in \cite[Section 12.5]{Nocedal_Wright} introduce the critical cone, containing directions orthogonal to the gradients of the constraint functions:
\[ \mathcal C(\bo x^*) = \{w \in \Bbb{R}^{2n} : \nabla \mathcal C_j(\bo x^*)\cdot w = 0,\ 0 \leq j \leq n-1\}.\]
The Lagrangian, ignoring constraints that are not active, is defined by $\mathcal L : \Bbb{R}^{2n}\times \Bbb{R}^n \to \Bbb{R}$ 
\[  \mathcal L(\bo x,\lambda) = |D(x_0,...,x_{n-1})|+\sum_{i=0}^{n-1} \lambda_i \mathcal C_i(\bo x).
\]  
Assume that  $\bo x^*$ is a local minimizer. Then $\bo x^*$ verifies the optimality condition \eqref{eq:lagrange-first} with corresponding Lagrange multipliers $\lambda^*= (\lambda_0,...,\lambda_{n-1})$ and the following second order condition holds
\[ w^T \nabla_{xx} \mathcal L(\bo x^*,\lambda^*) w \geq 0 \text{ for every } w \in \mathcal C(\bo x^*).\]
For more details see \cite[Theorem 12.5]{Nocedal_Wright}.

To compute the Hessian of the Lagrangian with respect to the $\bo x$ variable we proceed in a more direct approach. Observe that the disk-polygon $D(x_0,...,x_{n-1})$ is composed of the polygon $x_0...x_{n-1}$ and the unit circle segments given by chords of lengths $|x_i-x_{i+1}|$, $i=0,...,n-1$. This renders the area of a disk-polygon explicit in terms of the coordinates of the vertices.
  
The area of a segment of a unit disk corresponding corresponding to a central angle $\theta \in [0,\pi]$ is given by $\frac{\theta-\sin \theta}{2}$. Thus the area of a disk segment bounded by a chord of length $x$, corresponding to a central angle  $\theta = 2\arcsin \frac{x}{2}$ is
 \[ f(x):=\arcsin (x/2)-x/2\sqrt{1-x^2/4}.\]
 
 The first and second derivatives of $f$ are given by
 \[ f'(x) = \frac{x^2}{2\sqrt{4-x^2}}, \ \ f''(x) = \frac{x(8-x^2)}{2(4-x^2)^{3/2}}.\]
 In particular $f$ is convex on $[0,1]$.
 
Given two points $a\neq b \in \Bbb{R}^2$ and $f$ a real function of class $C^2$ consider the function $g: (a,b)\mapsto f(|a-b|)$. The first and second derivatives of this function with respect to $a,b$ are given by:
\begin{align}
& \partial_a g(a,b) = f'(|a-b|)\frac{a-b}{|a-b|},\ \partial_b g(a,b) = f'(|a-b|)\frac{b-a}{|a-b|}\notag \\
& \partial_{aa} g(a,b) =\partial_{bb} g(a,b)= f''(|a-b|)\frac{a-b}{|a-b|}\otimes \frac{a-b}{|a-b|}+\frac{f'(|a-b|)}{|a-b|}\left(I-\frac{a-b}{|a-b|}\otimes \frac{a-b}{|a-b|}\right)\notag \\
& \partial_{ab} g(a,b) = \partial_{ba}g(a,b) = -\partial_{aa}g(a,b)=-\partial_{bb}g(a,b).  \label{eq:hess-distance}
\end{align}

 With these considerations the Lagrangian becomes
 \[ \mathcal L(\bo x, \lambda) =  \text{Area}(x_0...x_{n-1})+\sum_{i=0}^{n-1} f(|x_i-x_{i+1}|)+\sum_{i=0}^{n-1}\lambda_i |x_i-x_{i+k}|. \]
 Let us compute the Hessian for all the terms involved for a regular Reuleaux $n$-gon. 
 
 \begin{itemize}
 	\item 
 The area of the polygon is given by the well known formula
 \[ \Area (\bo x) = \frac{1}{2} \sum_{i=0}^{n-1} (x_i^1x_{i+1}^2-x_{i+1}^1x_i^2).\]
 Thus, the Hessian matrix has the following block structure 
 \begin{equation}
 \nabla_{xx} \Area(\bo x) = \begin{pmatrix}
 \bo B_{ij}
 \end{pmatrix}_{0\leq i,j\leq n-1} 
 \label{eq:Hess_area}
 \end{equation}
 where the non-zero $2\times 2$ blocks are given by $\bo B_{ij} = \begin{pmatrix}0& \frac12 \\ -\frac12 & 0 \end{pmatrix}$ if $j = i+1$ and $\bo B_{ij} = \begin{pmatrix}0& -\frac12 \\ \frac12 & 0 \end{pmatrix}$ if $i = j+1$. 
 
 \item In view of the previous computations, the second term has the following contributions for $0 \leq i \leq n-1$:
 \begin{align*}\partial_{x_ix_i} f(|x_i-x_{i+1}|)& = \frac{|x_i-x_{i+1}|}{2(4-|x_i-x_{i+1}|^2)^{1/2}}I+\frac{2|x_i-x_{i+1}|}{(4-|x_i-x_{i+1}|^2)^{3/2}} \frac{x_i-x_{i+1}}{|x_i-x_{i+1}|}\otimes \frac{x_i-x_{i+1}}{|x_i-x_{i+1}|}\\
 & = \frac{1}{2}\tan \frac{\pi}{2n} I +\frac{\sin \frac{\pi}{2n}}{2\cos^3 \frac{\pi}{2n}} \frac{x_i-x_{i+1}}{|x_i-x_{i+1}|}\otimes \frac{x_i-x_{i+1}}{|x_i-x_{i+1}|}
 \end{align*}
 \[\partial_{x_{i+1}x_{i+1}} f(|x_i-x_{i+1}|) = \partial_{x_ix_i} f(|x_i-x_{i+1}|),  \]
 \[\partial_{x_{i}x_{i+1}} f(|x_i-x_{i+1}|) =\partial_{x_{i+1}x_{i}} f(|x_i-x_{i+1}|) = - \partial_{x_ix_i} f(|x_i-x_{i+1}|). \]
 \item The contributions of the third term are given by
 \begin{align*} \partial_{x_ix_i}\lambda_i |x_i-x_{i+k}| & = \lambda_i \frac{1}{|x_i-x_{i+k}|}\left(I- \frac{x_i-x_{i+k}}{|x_i-x_{i+k}|}\otimes\frac{x_i-x_{i+k}}{|x_i-x_{i+k}|}\right)\\
 & = -\tan \frac{\pi}{2n} \left(I- (x_i-x_{i+k})\otimes(x_i-x_{i+k})\right)
 \end{align*}
 \[ \partial_{x_{i+k}x_{i+k}}\lambda_i |x_i-x_{i+k}| = \partial_{x_ix_i}\lambda_i |x_i-x_{i+k}|, \]
 \[ \partial_{x_ix_{i+k}}\lambda_i |x_i-x_{i+k}| = \partial_{x_{i+k}x_{i}}\lambda_i |x_i-x_{i+k}| = -\partial_{x_ix_i}\lambda_i |x_i-x_{i+k}| \] 
  \end{itemize} 
 On the other hand, the points $x_i$ are the vertices of a Reuleaux triangle if and only if $|x_i-x_j|\leq 1$ for every $i$, and in particular $|x_i-x_j| = 1$ if and only if $j \in \{i+k,i-k\}$ modulo $n$. 
 
 Elements $\bo w=(w_i)_{i=0}^{n-1}$ belonging to the critical cone $\mathcal C(\bo x^*)$ verify
 \begin{equation}\label{eq:constr-w} (x_i-x_{i+k})\cdot(w_i-w_{i+k}) = 0\text{ for every } 0 \leq i \leq n-1.
 \end{equation}
Therefore, using the block structure of the Hessian of the Lagrangian found above we find: 
\begin{align*}
\bo w^T\nabla_{xx}\mathcal L(\bo x^*,\lambda) \bo w & = \sum_{i=0}^{n-1} w_i \begin{pmatrix}
0& 1\\-1 & 0
\end{pmatrix}w_{i+1} + \sum_{i=0}^{n-1} \frac{1}{2} \tan \frac{\pi}{2n} |w_i-w_{i+1}|^2+\\
& + \sum_{i=0}^{n-1} \frac{\sin \frac{\pi}{2n}}{2\cos^3 \frac{\pi}{2n}4\sin^2 \frac{\pi}{2n}} [(x_i-x_{i+1})\cdot(w_i-w_{i+1})]^2 \\
&-\tan \frac{\pi}{2n} \sum_{i=0}^{n-1} |w_i-w_{i+k}|^2.
\end{align*}

Denote $D_i = x_{i+k}-x_i$ (a unit vector in $\Bbb{R}^2$) and $W_i = w_{i+k}-w_i$. Denote by $v^\perp$ the rotation of the vector $v$ around the origin with angle $\pi/2$. Constraints \eqref{eq:constr-w} for perturbations belonging to the critical cone imply the existence of constants $q_0,...,q_{n-1} \in \Bbb{R}$ such that 
\[ W_i = q_i D_i^\perp.\]
It is straightforward to see that
\[ x_{i+1}-x_i = -D_{i+1}-D_{i-k}, w_{i+1}-w_i = -W_{i+1}-W_{i-k},\]
where indices are considered modulo $n$. We find that
\begin{align*}
|w_{i+1}-w_i|^2 & = q_{i+1}^2+q_{i-k}^2-2q_{i+1}q_{i-k} \cos \frac{\pi}{n} \\
\left|(x_{i+1}-x_i)\cdot (w_{i+1}-w_i)\right|^2 & = (-q_{i-k}\sin \frac{\pi}{n}+q_{i+1} \sin \frac{\pi}{n})^2=\sin^2 \frac{\pi}{n} (q_{i-k}-q_{i+1})^2
\end{align*}
With these considerations we find
\begin{align*}
\bo w^T\nabla_{xx}\mathcal L(\bo x^*,\lambda) \bo w & = \sum_{i=0}^{n-1} w_i \begin{pmatrix}
0& 1\\-1 & 0
\end{pmatrix}w_{i+1} + \frac{1}{2} \tan \frac{\pi}{2n}\sum_{i=0}^{n-1} (q_{i+1}^2+q_{i-k}^2-2q_{i+1}q_{i-k}\cos \frac{\pi}{n})\\
& +\frac{1}{2} \tan \frac{\pi}{2n} \sum_{i=0}^{n-1} (q_{i+1}-q_{i-k})^2 \\
&-\tan \frac{\pi}{2n} \sum_{i=0}^{n-1} q_i^2\\
& = \sum_{i=0}^{n-1} w_i \begin{pmatrix}
0& 1\\-1 & 0
\end{pmatrix}w_{i+1} + \frac{1}{2} \tan \frac{\pi}{2n}\left(2\sum_{i=0}^{n-1}q_i^2-2(\cos \frac{\pi}{n}+1)q_iq_{i+k}\right)
\end{align*}
Ideally one would like to prove that 
\[\bo w^T\nabla_{xx}\mathcal L(\bo x^*,\lambda) \bo w\leq 0\]
for every $\bo w \in \mathcal C(\bo x^*)$, showing that regular Reuleaux polygons are local maximizers for $n \geq 5$. However, the associated computations are not straightforward.

Nevertheless, for a Blaschke transformation (see Figure \ref{fig:perturbedR}) of the vertices $x_0,x_{-k}$ we have $q_0=q_1=1$ and $q_{-k} = 2\cos \frac{\pi}{n}$, all other terms being zero. The term 
\[\sum_{i=0}^{n-1} w_i \begin{pmatrix}
0& 1\\-1 & 0
\end{pmatrix}w_{i+1}
\]
cancels since only neighboring vertices for which $w_i\neq 0$ interact.

This leads to 
\[\bo w^T\nabla_{xx}\mathcal L(\bo x^*,\lambda) \bo w = 2\tan \frac{\pi}{2n}(1-2\cos \frac{\pi}{n})<0,\]
for $n \geq 5$. Therefore regular Reuleaux polygons are not local minima for $n \geq 5$, implying that the only possible minimizer is the Reuleaux triangle.

\subsection{Centers as variables} In the previous section the vertices of a disk polygon were taken as variables and the constant width constraint was imposed by requiring a family of distances between them to be equal to one. It is possible to take a different perspective, since a disk-polygon is completely characterized by its vertices, but also by the centers of the disks defining it through their intersection. See Figure \ref{fig:disk-polygon}.

Recall form Definition \ref{def:disk-poly} that if $c_0,...,c_{n-1}$ are points in $\Bbb{R}^2$ then $D = \bigcap_{i=0}^{n-1}B(c_i)$ is a disk-polygon. Unlike the case of Reuleaux polygons, the disk centers $c_i$ are not necessarily vertices of the disk-polygon. In the following we always assume that no disk is redundant in the definition of $D$.

A disk-polygon is a Reuleaux polygon if the constraints
\begin{equation}\label{eq:constr-centers}
\mathfrak C = \{(c_i)_{i=0}^{n-1} \in \Bbb{R}^{2n} : |c_i-c_{i+k}|=1, \text{ for }i=0,...,n-1; |c_i-c_j|\leq 1, \forall i\neq j\}
\end{equation} 
hold. Then the Blaschke-Lebesgue problem is equivalent to 
\begin{equation}\label{eq:BL-centers}
 \min_{(c_i)\in \mathfrak C} \text{Area}(\bigcap_{i=0}^{n-1}B(c_i)).
 \end{equation}

Compared to the previous section, the variation of the area is easily expressed using the variation of the centers. The boundary of $\bigcap_{i=0}^{n-1}B(c_i)$ is made of arcs of circles with centers $c_i$, therefore, the infinitesimal normal movement is strictly given by the perturbations of the vertices $c_i$.

In the following, fix $n=2k+1$ and assume that $\bo c^* = (c_0,...,c_{n-1})$ is a solution for \eqref{eq:BL-centers}. Without loss of generality, assume that all centers are non-redundant. If this is not the case, eliminate the redundant centers and decrease accordingly the value of $n$. 
 
For each center $c_i$ denote by $\gamma_i$ the (non-void) arc of unit circle centered at $c_i$ contained in the boundary of the disk-polygon. The derivative of the area of the disk-polygon with respect to the variation the disk center $c_i$ is given by
\[ \frac{\partial \text{Area}(\bigcap_{i=0}^{n-1}B(c_i))}{\partial c_i} = \int_{\gamma_i} \bo n.\]
This is a direct consequence of the shape derivative of the area \eqref{eq:area-deriv}. Assuming a minimizer for \eqref{eq:BL-centers} exists and $\mathfrak C$ is not trivial (which is the case when $n \geq 5$) the classical optimality conditions state that the gradient of the objective function is a linear combination of the gradients of the active constraints.

Define the Lagrangian $\mathcal L : \Bbb R^{2n} \times \Bbb R^n \to \Bbb{R}$:
\[\mathcal L(x,\lambda) = \text{Area}(\bigcap_{i=0}^{n-1}B(c_i))+\sum_{i=0}^{n-1} \lambda_i |c_i-c_{i+k}|.\]

\bo{First order optimality conditions.} For a minimizer the derivative with respect to $c_i$ of the Lagrangian vanishes giving:
\[ \frac{\partial \text{Area}(\bigcap_{i=0}^{n-1}B(c_i))}{\partial c_i} = \int_{\gamma_i} \bo n =-\lambda_i (c_i-c_{i+k})-\lambda_{i-k}(c_i-c_{i-k}). \]
 
Since the average normal on $\gamma_i$ is aligned with the symmetry axis and $(c_i-c_{i+k}),(c_i-c_{i-k})$ have unit length we immediately find that $\lambda_i=\lambda_{i-k}$ for every $i=0,...,n-1$, therefore all Lagrange multipliers are equal at a critical point. More precisely:
\[\int_{\gamma_i} \bo n = \tan \frac{\theta_i}{2} (c_{i+k}-c_i+c_{i-k}-c_i), \]
showing that $\lambda_i = \tan \frac{\pi}{n}$, for every $i=0,...,n-1$.

\begin{rem}
	Observe the following difference in the sign of the Lagrange multipliers for the two situations where vertices or centers are taken as variables:
	\begin{itemize}[noitemsep]
		\item \bo{Vertices as variables}: a critical point for the Lagrangian gives multipliers $\lambda_i= -\tan \frac{\pi}{2n}$. The constant width constraint keeps the vertices of a disk-polygon \emph{together} when minimizing the area. 
		\item \bo{Centers as variables}: a critical point for the Lagrangian gives multipliers $\lambda_i= \tan \frac{\pi}{2n}$. The constant width constraint keeps the disk centers for a disk-polygon \emph{apart} when minimizing the area. 
	\end{itemize}
\end{rem}

\bo{Second order optimality conditions.} For computations of the second derivative of the area of a disk-polygon with respect to center movement, we refer to the computations in \cite{Laurain2022}, formulas (2.8), (2.9). In the cited reference a more complex situation is considered, where the boundary of a domain consists of general arc of circles (not necessarily convex, different radii, etc.). Nevertheless, the formulas apply to the case of disk-polygons and plugging all the required information into the formulas from \cite{Laurain2022} gives the following $2\times 2$ block structure for the Hessian of the area of a disk-polygon with respect to the centers:
\begin{align*}
\partial_{c_ic_i}A(\bo c) & = \int_{\gamma_i} \bo n \otimes \bo n - \bo t \otimes \bo t  - \frac{1}{\tan \theta_{i-k}} (c_{i-k}-c_i)\otimes (c_{i-k}-c_i) \\
 & - \frac{1}{\tan \theta_{i+k}} (c_{i+k}-c_i)\otimes (c_{i+k}-c_i) \\
\partial_{c_ic_{i+1}}A(\bo c) & = \frac{1}{\sin \theta_{i-k}} (c_{i-k}-c_i)\otimes (c_{i-k}-c_{i+1})\\
\partial_{c_ic_{i-1}}A(\bo c) & = \frac{1}{\sin \theta_{i+k}} (c_{i+k}-c_i)\otimes (c_{i+k}-c_{i-1})
\end{align*}
where $\bo a \otimes \bo b = (a_ib_j)_{1\leq i,j \leq 2}$ is the classical tensor product and $A(\bo c) =\text{Area}(\bigcap_{i=0}^{n-1}B(c_i))$. All other blocks are zero in the Hessian matrix. Let us recall the following property which will justify some of the computations below: $\bo x^T (\bo a \otimes \bo b) \bo y = (\bo x \cdot \bo a)(\bo b\cdot \bo y)$.

Consider the Lagrangian 
\[ \mathcal L(\bo c, \lambda) = A(\bo c)+\sum_{i=0}^n \lambda_i|c_{i+k}-c_i|.\]
For $\bo c^*$ solution in \eqref{eq:BL-centers}, consider the associated  critical point for the Lagrangian $(\bo c^*,\lambda)$. Furthermore, consider perturbations $\bo w$ in the critical cone, verifying
\[ (c_{i+k}-c_i)\cdot (w_{i+k}-w_i) = 0,\ i=0,...,n-1.\]
Remembering the computations for Hessians of distance functions in \eqref{eq:hess-distance} gives 
\begin{align*} \bo w^T \partial_{\bo c\bo c}\mathcal L(\bo c,\lambda) \bo w & =  \bo w^T\partial_{\bo c\bo c}A(\bo c)\bo w + \sum_{i=0}^{n-1} \lambda_i|w_{i+k}-w_i|^2.\\
& = w_0\partial_{ c_0 c_0}A(\bo c)w_0+ w_{-k}\partial_{ c_{-k} c_{-k}}A(\bo c)w_{-k}+\lambda_0 |w_0|^2+\lambda_1 |w_{-k}|^2+\lambda_{-k}|w_0-w_{-k}|^2.
\end{align*}

Like in the previous case, to simplify the computations, we focus on Blaschke perturbations for centers $c_0$, $c_{-k}$, which have two advantages. First, all remaining $n-2$ centers are not perturbed and secondly, no neighboring vertices interact.

Some preliminary computations show that if $w_i$ is a perturbation for $c_i$ then:
\begin{align*}
w_i^T \left(\int_{\gamma_i}\bo n \otimes \bo n-\bo t \otimes \bo t\right) w_i &=\int_{\gamma_i} (w\cdot \bo n)^2-(w\cdot \bo t)^2\\
 & = \sin(\theta_i)((w_i\cdot \bo b_i)^2-(w\cdot \bo t_i)^2)
\end{align*}
where $\bo b_i$ is the unit bisector of $\angle c_{i-k}c_ic_{i+k}$ and $\bo t_i = \bo b_i^\perp$.

Consider {Blaschke perturbations} for centers $c_0, c_{-k}$ with perturbations $w_0 \perp (c_0-c_k)$, $w_{-k}\perp (c_1-c_{-k})$, all other $w_i, i\notin \{0,-k\}$ being zero. Since $(w_0-w_{-k})\cdot (c_0-c_{-k})=0$ we find $|w_0|\sin \theta_0 = |w_{-k}|\sin \theta_{-k}=K$. Also, geometrically it is immediate that $\angle w_0,w_{-k} = \pi-\frac{\theta_0+\theta_{-k}}{2}$.

Therefore
\begin{align*}
\bo w \partial_{cc}\mathcal L(\bo c^*,\lambda) \bo w & =-\frac{2K^2}{\tan \theta_0}-\frac{2K^2}{\tan \theta_{-k}} + \lambda_0 |w_0|^2+\lambda_1|w_{-k}|^2+\lambda_{-k}|w_{-k}-w_0|^2\\ 
& = K^2\Big(-\frac{2}{\tan \theta_0}-\frac{2}{\tan \theta_{-k}}+\frac{\lambda_0}{\sin^2 \theta_0}+\frac{\lambda_1}{\sin^2 \theta_{-k}}\\
&+2\lambda_{-k}\cos \frac{\theta_0+\theta_{-k}}{2} \frac{1}{\sin \theta_0 \sin \theta_{-k}}\Big).
\end{align*}

For a critical shape for the Lagrangian we have $\lambda_i = \tan \frac{\pi}{2n}$ and $\theta_i = \pi/n$ for $i=0,...,n-1$.
Replacing these quantities in the formula above gives
\[ \bo w \partial_{cc}\mathcal L (\bo c^*,\lambda) \bo w \sim (1-2\cos \frac{\pi}{n}),\]
where $\sim$ is understood up to a multiplication with a positive constant. Therefore, for $n \geq 5$ any regular Reuleaux $n$-gon is not a local minimum since along a Blaschke perturbation we obtain  $\bo w \partial_{cc}\mathcal L(\bo c^*,\lambda)\bo w<0$.

\section{Conclusions}

Two purely variational arguments (direct and Lagrangian) used in this paper to show that the only critical Reuleaux polygons for the area functional are the regular ones. Moreover, regular Reuleaux $n$-gons are not local minimizers, since there exist perturbations which decrease the area. Thus, the Reuleaux triangle minimizes the area among all shapes having constant width. This is in accord with previous observations in \cite{Harrell}, \cite{henrot_lucardesi_annulus}.

It can also be observed that the minimality of the Reuleaux triangle is a consequence of its extremal character: none of its vertices can be perturbed like in Theorem \ref{thm:sensitivity}. It is the only \emph{rigid} Reuleaux $n$-gon. Moreover, no local minimizers for the area exist among Reuleaux $n$-gons with $n \geq 5$. This shows why many proof strategies for the Blaschke-Lebesgue theorem rely simply on modifying a given Reuleaux $n$-gon to obtain another one with a smaller area, which is always possible if $n \geq 5$.

Ideas presented in this paper could lead to new developments regarding the three dimensional problem, still open today \cite{kawohl-webe}: 

\emph{Which three dimensional body minimizes the volume among bodies having constant width?}

The three dimensional analogue of the Reuleaux polygons, the Meissner polyhedra, were introduced recently in \cite{montejano} and analyzed further in \cite{meissner_hynd}. The volume of these Meissner polyhedra is computed explicitly in \cite{bogosel_Meissner} and \cite{hynd-vol-per}. The missing ingredient for mimicking the arguments presented in this paper in the three dimensional case is understanding perturbations of Meissner polyhedra. The main difficulty is the fact that the combinatorial structure of the diametrical pairs is more complex than in dimension two (where each vertex is linked with exactly two diameters). The Lagrangian perspective presented in Section \ref{sec:Lagrangian} could help in this direction, since the combinatorial difficulties can be simplified allowing independent movement of opposite endpoints of a diameter.

\bibliographystyle{abbrv}
\bibliography{./biblio}

\bigskip
\small\noindent
Beniamin \textsc{Bogosel}: Centre de Math\'ematiques Appliqu\'ees, CNRS,\\
\'Ecole Polytechnique, Institut Polytechnique de Paris,\\
91120 Palaiseau, France \\
{\tt beniamin.bogosel@polytechnique.edu}\\
{\tt \nolinkurl{http://www.cmap.polytechnique.fr/~beniamin.bogosel/}}

\end{document}